\documentclass[12pt]{article}
\usepackage{amsthm, amssymb}
\newcommand{\nnn}[1]{(\ref{#1})}
\newcommand{\operatorname}[1]{{\rm #1}}
\newcommand{\dfrac}[2]{\displaystyle{\frac{#1}{#2}}}
\newcommand{\tfrac}[2]{\textstyle{\frac{#1}{#2}}}

\newcommand{\twostac}[2]{\left(\begin{array}{c}{#1}\\{#2}\end{array}\right)}
\newcommand{\HTmx}[2]{\left(\begin{array}{c}{#1}\\{#2}\end{array}\right)}
\newcommand{\htmx}[2]{\begin{array}{c}{#1}\\{#2}\end{array}}

\newcommand{\HHmx}[4]{\left(\begin{array}{cc}{#1}&{#2}\\{#3}&{#4}\end{array}
\right)}
\newcommand{\HTTmx}[3]{\left(\begin{array}{c}{#1}\\{#2}\\{#3}\end{array}\right)}
\newcommand{\THmx}[2]{\left(\begin{array}{cc}{#1}&{#2}\end{array}\right)}

\newtheorem{lemma}{Lemma}[section]
\newtheorem{rmk}{Remark}

\let\a=\alpha
\let\b=\beta
\let\c=\chi
\let\d=\delta

\let\e=\varepsilon

\let\f=\varphi
\let\F=\Phi
\let\g=\gamma
\let\G=\Gamma

\let\i=\iota
\let\l=\lambda
\let\L=\Lambda
\let\m=\mu

\let\N=\nabla
\let\nd=\nabla

\let\r=\rho

\let\t=\tau

\let\vp=\varpi

\let\w=\omega

\let\X=\Xi

\let\y=\psi
\let\Y=\Psi

\def\bbN{{\mathbb N}}
\def\bbR{{\mathbb R}}

\def\bbV{{\mathbb V}}

\def\cD{{\cal D}}

\def\cL{{\cal L}}

\def\cP{{\cal P}}
\def\cR{{\cal R}}

\def\cV{{\cal V}}

\def\eye{\sqrt{-1}}

\def\spin{\operatorname{Spin}}

\def\tf{\tfrac12}

\def\tt{\tilde\theta}
\def\tta{\tilde\tau}

\let\ptl=\partial
\def\pt{\partial_t}

\def\sideremark#1{\ifvmode\leavevmode\fi\vadjust{\vbox to0pt{\vss% the remark
 \hbox to 0pt{\hskip\hsize\hskip1em%                          will appear only
 \vbox{\hsize2cm\tiny\raggedright\pretolerance10000%          on the side
 \noindent #1\hfill}\hss}\vbox to8pt{\vfil}\vss}}}%
\title{Spectrum Generating on Twistor Bundle}
\author{Thomas Branson and Doojin Hong}

\begin{document}

\maketitle

\abstract{We give explicit formulas for 
the intertwinors of all orders on the 
twistor bundle over $S^1\times S^{n-1}$ 
using spectrum generating technique introduced in \cite{BOO:96}.}

%%%%%%%%%%%%%%%%%%%%%%%%%%%%%%%%%%%%%%%%%%%%%%%%%%%%%%%%%%%%%%%
\section{Introduction}
%%%%%%%%%%%%%%%%%%%%%%%%%%%%%%%%%%%%%%%%%%%%%%%%%%%%%%%%%%%%%%%
It was shown in \cite{BOO:96} that one can construct intertwining 
operators of some representations without too much effort when 
eigenspaces occur with multiplicity one. On the differential 
form bundle over $S^1\times S^{n-1}$, the double cover of the compactified  
Minkowski space, some $K$-type eigenspaces occur with multiplicity two. 
After some additional computation, Branson also showed spectral function for these 
operators. 

Intertwinors on spinors like the Dirac operator have eigenspaces with 
multiplicity one over $S^1\times S^{n-1}$and explicit spectral function 
was given in \cite{Hong:04}. But on twistors, the eigenspaces of the 
intertwinors including Rarita Schwinger operator have multiplicity two 
on some $K$-type. In this paper, we present the spectral function for these 
operators.

We briefly review conformal covariance and intertwining 
relation (for more details, see 
\cite{Branson:96}, \cite{BOO:96}).
\\
Let $M$ be an n-dimensional spin manifold. We enlarge the 
structure group 
$\spin(n)$ to $\spin(n)\times \bbR_+$ in conformal geometry. 
$(V(\l),\l^r)$ are finite dimensional 
$\spin(n)\times \bbR_+$ representations, 
where $(V(\l),\l)$ are finite dimensional 
representations of $\spin(n)$ and 
$\l^r(h,\a)=\a^r\l(h)$ for $h\in \spin(n)$ 
and $\a\in \bbR_+$. The corresponding associated 
vector bundles are $\bbV(\l)=P_{\spin(n)}\times_{\l}V(\l)$ 
and $\bbV^r(\l)=P_{\spin(n)\times \bbR_+}\times_{\l^r}V(\l)$ 
with structure groups $\spin(n)$ and $\spin(n)\times \bbR_+$. $r$ 
is called the conformal weight of $\bbV^r$. 
Tangent bundle $TM$ carries conformal weight $-1$ 
and cotangent bundle $T^*M$ carries conformal weight $+1$. 
In general, if $V$ is a subbundle of 
$(TM)^{\otimes p}\otimes (T^*M)^{\otimes q}
\otimes(\Sigma M)^{\otimes r}\otimes 
(\Sigma^* M)^{\otimes s}$, then $V$ carries conformal weight 
$q-p$, where $\Sigma M$ is the contravariant spinor bundle. 
\\
A conformal covariant of bidegree $(a,b)$ 
is a $\spin(n)\times \bbR_+$-equivariant differential operator 
$D: \bbV^r(\l) \rightarrow \bbV^s(\sigma)$ 
which is a polynomial in the metric $g$, its inverse $g^{-1}$, 
the volume element $E$, and the 
fundamental tensor-spinor $\g$ with a conformal covariance law 
$$
\omega \in C^{\infty}, \quad \overline{g}=e^{2\omega}g, 
\quad \overline{E}=e^{n\omega}E, \quad \overline{\g}=e^{-\omega}\g 
\Rightarrow \overline{D}=e^{-b\omega}D\m(e^{a\omega}) ,
$$
where $\m(e^{a\omega})$ is multiplication of $e^{a\omega}$. 

\noindent Given a conformal covariant of bidegree 
$(a,b)$, $D:\bbV^r(\l) \rightarrow \bbV^s(\sigma)$, 
we can assign new conformal 
weights to get $D:\bbV^{r'}(\l) \rightarrow \bbV^{s'}(\sigma)$ 
whose bidegree is then $(a-r'+r,b-s'+s)$. Calling this 
$D$ again is an abuse of notation. 
If $r'=r+a$ and $s'=s+b$, then $D:\bbV^{r+a}(\l) 
\rightarrow \bbV^{s+b}(\sigma)$ becomes conformally invariant and we 
call $(a+r,b+r)$ the reduced conformal bidegree of $D$. 
To see how conformal covariants behave under a conformal transformation
and a conformal vector field, we recall followings.\\
A diffeomorphism $h:M\rightarrow M$ is called a 
conformal transformation if $h\cdot g=e^{2\omega_h}g$, where $\cdot$ is the 
natural action of $h$ on tensor fields.
A conformal vector field is a vector field $X$ with $\cL_Xg=2\omega_Xg$ for 
some $\omega_X \in C^{\infty}(M)$. 
A conformal covariant $D:\bbV^0(\l) \rightarrow \bbV^0(\sigma)$ 
of reduced bidegree 
$(a,b)$ satisfies
$$
D(e^{a\omega_h}h\cdot\f)=e^{b\omega_h}h\cdot (D(\f)) 
\quad {\rm and} \quad D(\cL_X+a\omega_X)\f=(\cL_X+b\omega_X)D\f .
$$
for all $\f \in \Gamma(\bbV^0(\l))$. Thus if 
$D:\bbV^r(\l) \rightarrow \bbV^s(\sigma)$ of reduced bidegree 
$(a,b)$, then 
\begin{equation}\label{inter}
D(\cL_X+(a-r)\omega_X)\f=(\cL_X+(b-s)\omega_X)D\f
\end{equation}
for $\f\in \Gamma(\bbV^r(\l))$ and $D\f\in 
\Gamma(\bbV^s(\sigma))$.\\
Note that conformal vector fields form a Lie algebra 
$\mathfrak{c}(M,g)$ and give rise to the 
principal series representation 
$$
U_a^{\l}:\mathfrak{c}(M,g)\rightarrow 
{\rm End}\Gamma(\bbV^0(\l)) \quad {\rm by} \quad X \mapsto \cL_X+a\omega_X .
$$
So a conformal covariant $D:\bbV^r(\l) 
\rightarrow \bbV^s(\sigma)$ of reduced bidegree 
$(a,b)$ intertwines the principal series representation
$$
DU_{a-r}^{\l}\f=U_{b-s}^{\sigma}D\f
$$
for $\f\in \Gamma(\bbV^r(\l))$ and $D\f\in \Gamma(\bbV^s(\sigma))$.

%%%%%%%%%%%%%%%%%%%%%%%%%%%%%%%%%%%%%%%%%%%%%%%%%%%%%%%%%%%%%%%%%%%%
\section{Spinors and Twistors}
%%%%%%%%%%%%%%%%%%%%%%%%%%%%%%%%%%%%%%%%%%%%%%%%%%%%%%%%%%%%%%%%%%%%
Let $M=S^1\times S^{n-1}$, $n$ even, be a manifold endowed with the 
Lorentz metric $-dt^2+g_{S^{n-1}}$.

To get a fundamental tensor-spinor $\a$ for $M$ from the corresponding
object $\g$ on $S^{n-1}$, let
$$
\a^j=\left(\begin{array}{cc}
\g^j & 0 \\
0 & -\g^j
\end{array}\right)\,,\qquad j=1,\ldots,n-1,
$$
and 
$$
\a^0=\left(\begin{array}{cc}
0 & 1 \\
1 & 0
\end{array}\right)\,.
$$

Since $M$ is even-dimensional, there is a {\em chirality} operator $\c_M\,$,
equal to some complex unit times $\a^0\tilde\c_S\,$, where 
$$
\tilde\c_S=\left(\begin{array}{cc}
\c_S & 0 \\
0 & -\c_S
\end{array}\right)\,,
$$
$\c_S$ being the  
chirality operator on $S$.  The chirality operator is always normalized
to have square $1$; thus $(\c_S)^2$ and
$(\tilde\c_S)^2$ are identity operators, and since $\a^0\a^0=1$,
we have
$(\a^0\tilde\c_S)^2=-1$.
As a result, we may take
$$
\c_M=\pm\sqrt{-1}\a^0\tilde\c_S.
$$
A {\em spinor} on $M$ can be viewed as a pair of time-dependent spinors 
on $S^{n-1},$ i.e.,
$
\left(\begin{array}{c}
	   \varphi \\
	   \psi
      \end{array} 
\right),
$
where $\varphi$ and $\psi$ are $t$-dependent spinors on $S^{n-1}.$ 
But by chirality consideration (\cite{Hong:00}), 
we get $\Xi=\pm 1$ spinors:
$$
\left(\begin{array}{c}
	   \varphi \\
	   \psi
      \end{array} 
\right)
=
\left(\begin{array}{c}
	   \Xi\psi/\eye \\
	   \psi
      \end{array} 
\right)\, .
$$
Recall that {\em twistors} are spinor-one-forms $\F_\l$ with $\a^\l\F_\l=0$.
Given a chirality $\Xi$, a twistor $\Y$ is determined by a $t$-dependent 
spinor-one-form $\y_j$ on $S^{n-1}$ via
$$
\Y=dt\wedge\twostac{\f_0}{\y_0}+\twostac{\f_j}{\y_j}\,,
$$
where
$$
\begin{array}{rl}
\f_j&=-\Xi\sqrt{-1}\y_j, \\
\y_0&=\Xi\sqrt{-1}\g^k\y_k, \\
\f_0&=\g^k\y_k\,.
\end{array}
$$
Furthermore, by Hodge theoretic consideration (\cite{Hong:00}), 
twistors on $M$ can be decomposed into three pieces
\begin{equation}\label{decom}
\left( \begin{array}{cc}
             -(n-1) \theta & -\X\eye\g_{i}\theta \\
	     -(n-1)\X\eye\theta & \g_{i}\theta 
       \end{array} \right) 
+ \left( \begin{array}{cc}
	            0 & -\X\eye T_{i}\tau \\
	            0 &  T_{i}\tau 
	           \end{array} \right)
+   \left( \begin{array}{cc}
	            0 & -\X\eye\nabla^{j}\eta_{ji} \\
	            0 & \nabla^{j}\eta_{ji}
	           \end{array} \right) 
\end{equation}
$$
=:\langle\theta\rangle+\{\tau\}+[\eta]\, .
$$
%%%%%%%%%%%%%%%%%%%%%%%%%%%%%%%%%%%%%%%%%%%%%%%%%%%%%%%%%%%%%%%%%%%%%%%%%%%%%
\section{Intertwining relation on twistors}
%%%%%%%%%%%%%%%%%%%%%%%%%%%%%%%%%%%%%%%%%%%%%%%%%%%%%%%%%%%%%%%%%%%%%%%%%%%%%
Consider the standard conformal vector field 
(\cite{Branson:87, Orsted:81}) 
$$
T:=\cos\r\sin t\partial_t+\cos t\sin\r\partial_\r\, .
$$  
Here $\r$ is the azimuthal angle
on $S^{n-1}$.  The conformal factor of $T$ is
$$
\varpi:=\cos t\cos\r.
$$
Let $A=A_{2r}$ be an intertwinor of order $2r$. 
The intertwining relation says ((\ref{inter}), \cite{Branson:96, 
Branson:97, BOO:96}) 
\begin{equation}\label{rel-1}
A\left(\tilde\cL_T+\left(\frac{n}2-r\right)\varpi\right)=
\left(\tilde\cL_T+\left(\frac{n}2+r\right)\varpi\right)A,
\end{equation}
where $\tilde\cL_T$ is the {\em reduced Lie derivative}.
On a tensor-spinor with $\twostac{p}{q}$ tensor content, 
this is
$$
\tilde\cL_T=\cL_T+(p-q)\varpi.
$$
So here (with only 1-form content), it is $\cL_T-\varpi$.  
Note that we are
using the convention where spinors do not have an internal 
weight; otherwise
the spinor content would influence the reduction.\\
Since intertwinors change chirality, we want to consider  
an exchange operator
$$
\begin{array}{rl}
E:&=\a^0(\i(\partial_t)\e(dt)-\e(dt)\i(\partial_t)) \\
&=\a^0(1-2\e(dt)\i(\partial_t)).
\end{array}
$$
It is immediate that $E^2={\rm Id}$.
Because of the $\a^0$ factor, $E$ reverses chirality.
To see that $E$ takes twistors to twistors, note that
$$
\i(\partial_t)\e(dt)
-\e(dt)\i(\partial_t):\Phi_\l\mapsto\Phi_\l
-2\delta_\l{}^0\Phi_0\,.
$$
Thus 
$$
\begin{array}{rl}
\a^\l(E\F)_\l&=\a^\l\a^0(\F_\l-2\d_\l{}^0\F_0) \\
&=-2g^{\l 0}(\F_\l-2\d_\l{}^0\F_0)+2\a^0\a^\l\d_\l{}^0\F_0 \\
&=\underbrace{-2\F^0}_{2\F_0}
+4\underbrace{g^{00}}_{-1}\F_0+2\underbrace{\a^0\a^0}_1\F_0 \\
&=0,
\end{array}
$$
as desired.\\
We want to convert the relation (\ref{rel-1}) for $EA$.  
So we will eventually need $\cL_TE$.  We have:
$$
\begin{array}{rl}
\cL_TE&=\cL_T\left\{\a(dt)(1-2\e(dt)\i(\partial_t))\right\} \\
&=\{-\varpi\a(dt)+\a(d(Tt))\}(1-2\e^0\i_0) \\
&\qquad-2\a^0\{\e(dt)\i([T,\partial_t])+\e(d(Tt)\i(\partial_t)\}.
\end{array}
$$
But
$$
\begin{array}{rl}
Tt&=\cos\r\sin t, \\
d(Tt)&=-\sin\r\sin t\,d\r+\cos\r\cos t\,dt, \\
{}[T,\partial_t]&=-\cos\r\cos t\,\partial_t+\sin t\sin\r\,\partial_\r\,.
\end{array}
$$
This reduces the above to
\begin{equation}\label{lte}
\begin{array}{rl}
\cL_TE&=\sin t\a(d\w)(1-2\e^0\i_0)-2\sin t\a^0(\e^0\i(Y)+\e(d\w)\i_0) \\
&=\sin t\sin\r\{-\a^1(1-2\e^0\i_0)-2\a^0(\e^0\i_1-\e^1\i_0)\}.
\end{array}
\end{equation}
By Kosmann (\cite{Kosmann:72}, eq(16)), the Lie and covariant 
derivatives on spinors are related by
$$
\cL_X-\nd_X=-\tfrac14\nd_{[a}X_{b]}\g^a\g^b
=-\tfrac18(dX_\flat)_{ab}\g^a\g^b.
$$
Note that 
$$
\begin{array}{rl}
T_\flat&=-\cos\r\sin t\,dt+\cos t\sin\r\,d\r, \\
dT_\flat&=2\sin\r\sin t\,d\r\wedge dt.
\end{array}
$$
and
$$
d\varpi=-T_{\flat,{\rm R}}\,,
$$
where $\flat,$R is the musical isomorphism in the ``Riemannian" metric.
According to the above, 
\begin{equation}\label{ltandnt}
\cL_T-\nd_T=-\frac12\sin\r\sin t\a^1\a^0
\end{equation}
on spinors.\\
On a 1-form $\eta$,
$$
\langle(\cL_T-\nd_T)\eta,X\rangle=-\langle\eta,(\cL_T-\nd_T)X\rangle,
$$
since $\cL_T-\nd_T$ kills scalar functions.
But by the symmetry of the pseudo-Riemannian connection,
$$
[T,X]-\nd_TX=-\nd_XT.
$$
We conclude that 
$$
(\cL_T-\nd_T)\eta=\langle\eta,\nd T\rangle,
$$
where in the last expression, $\langle\cdot,\cdot\rangle$ is the pairing
of a 1-form with the contravariant part of a $\twostac{1}{1}$-tensor:
$$
((\cL_T-\nd_T)\eta)_\l=\eta_\m\nd_\l T^\m.
$$
Combining this with what we derived above for 
spinors (\ref{ltandnt}), 
for a spinor-1-form
$\F_\l$, we have
$$
((\cL_T-\nd_T)\F)_\l=\F_\m\nd_\l T^\m
-\frac12\sin\r\sin t\a^1\a^0\F_\l\,.
$$
But $\nd T$ {\em a priori} has projections in 3 
irreducible bundles, TFS${}^2$,
$\L^0$, and $\L^2$ (after using the musical isomorphisms).  
By conformality,
the TFS${}^2$ part is gone.  We expect a $\L^0$ part, 
essentially $\varpi$.
We also found the $\L^2$ part above,
$$
dT_\flat=2\sin\r\sin t\,d\r\wedge dt.
$$
More precisely, tracking the normalizations,
$$
(\nd T_\flat)_{\l\m}=(\nd T_\flat)_{(\l\m)}
+(\nd T_\flat)_{[\l\m]}
=(\varpi g+\frac12dT_\flat)_{\l\m}\,.
$$
Now note that 
$$
\begin{array}{rl}
\F_\m\nd_\l T^\m&=(\nd_\bullet T^\bullet\sharp\F)_\l \\
&=\varpi(g\sharp \F)_\l
+\frac12((dT_\flat)_{\nu\m}\e^\nu\i^\m\F)_\l \\
&=\varpi\F_\l+\frac12(((dT_\flat)_{01}\e^0\i^1
+(dT_\flat)_{10}\e^1\i^0)\F)_\l \\
&=\varpi\F_\l+\frac12((-2\sin\r\sin t\e^0\i^1
+2\sin\r\sin t\e^1\i^0)\F)_\l \\
&=\varpi\F_\l-\sin\r\sin t((\e^0\i^1-\e^1\i^0)\F)_\l \\
&=\varpi\F_\l-\sin\r\sin t((\e^0\i_1+\e^1\i_0)\F)_\l\,.
\end{array}
$$
As a result,
$$
\begin{array}{rl}
\cL_T-\nd_T&=\varpi-\sin\r\sin t\left(\tfrac12\a^1\a^0+\e^0\i_1
+\e^1\i_0\right) \\
&=:\varpi-\sin\r\sin t P \\
&=:\varpi-\cP,
\end{array}
$$
and 
$$
\tilde\cL_T-\nd_T=-\cP.
$$
An explicit calculation using \nnn{lte} gives
$$
(\cL_TE)E=-2\cP.
$$
Since $E^2={\rm Id}$, we conclude that 
$$
\cL_TE=-2\cP E.
$$
With the above, the intertwining relation for $EA$ becomes
$$
\begin{array}{rl}
\left(\tilde\cL_T+\left(\frac{n}2+r\right)\varpi\right)EA
&=E\left(\tilde\cL_T+\left(\frac{n}2+r\right)\varpi\right)A+(\cL_TE)A \\
&=EA\left(\tilde\cL_T+\left(\frac{n}2-r\right)\varpi\right)-2\cP EA,
\end{array}
$$
so that, with $B=EA$,
$$
B\left(\nd_T+\left(\frac{n}2-r\right)\varpi-\cP\right)=
\left(\nd_T+\left(\frac{n}2+r\right)\varpi+\cP\right)B.
$$
To see what $P$ does let us define two convenient operations.
$$
\y_j\stackrel{{\bf expa}}{\longmapsto}
\HHmx{u}{\X\y_j/\eye}{-\X u/\eye}{\y_j}
\stackrel{{\bf slot}}{\longmapsto}{\y_j}\, ,
$$
where $u=\g^k\y_k$.\\
Note that
$$
\begin{array}{l}
\y_j\stackrel{{\bf expa}}{\longmapsto}
\HHmx{u}{\X\y_j/\eye}{-\X u/\eye}{\y_j}
\stackrel{\i_0}{\longmapsto}\HTmx{u}{-\X u/\eye} \\
\stackrel{\e^1}{\longmapsto}
\HHmx{0}{\e^1 u}{0}{-\X\e^1 u/\eye}
\stackrel{{\bf slot}}\longmapsto-\X\e^1 u/\eye.
\end{array}
$$
As for the $\e^0\i_1$ term, anything in the range of $\e^0$ has a 
{\bf slot} of $0$.

Finally,
$$
\begin{array}{l}
\y_j\stackrel{{\bf expa}}{\longmapsto}
\HHmx{u}{\X\y_j/\eye}{-\X u/\eye}{\y_j}
\stackrel{\a^0}{\longmapsto}\HHmx{0}{1}{1}{0}
\HHmx{u}{\X\y_j/\eye}{-\X u/\eye}{\y_j} \\
=\HHmx{-\X u/\eye}{\y_j}{u}{\X\y_j/\eye}
\stackrel{\a^1}{\longmapsto}
\HHmx{-\X\g^1 u/\eye}{\g^1\y_j}{-\g^1 u}{-\X\g^1\y_j/\eye}.
\end{array}
$$
So
$$
\begin{array}{l}
{\bf slot}\,P\,{\bf expa}:\y_j\mapsto
-\frac12\X\g^1\y_j/\eye-\X(\e^1 u)_j/\eye \\
=-\frac{\X}{\eye}(\frac12\g^1\y_j+(\e^1 u)_j)
=-\frac{\X}{\eye}(\frac12\g^1\y_j+\d_j{}^1 u).
\end{array}
$$

Up to a factor of a complex unit, ${\bf slot}\,P\,{\bf expa}$ is
$$
\frac12\g^1\y_j+\d_j{}^1\g^k\y_k\,.
$$
We can also get this expression by successively taking the commutator of 
$\varpi$ with $\pt$ and
$$
{\bf slot}\,\cD\,{\bf expa}:\y_j\mapsto\frac12\g^k\N_k\y_j+\g^k\N_j\y_k\,.
$$
That is,
$$
\cP=\Xi\eye[\pt,[\cD,\varpi]]\,.
$$
Recall that $\cP=\sin \r \sin t P$.\\
After some straightforward computation, we get the block matrix for 
$\cD$ relative to the decomposition  
$\{\langle \theta \rangle, \{\tau \}, [\eta] \}$ (\ref{decom}) 
as follows.
$$
\left(\begin{array}{ccc}
	\displaystyle{\frac{n+1}{2(n-1)}J_\theta} & 
	\displaystyle{\frac{n-2}{4}-\frac{n-2}{(n-1)^2}J_\t^2} & 0\\
	-n & \displaystyle{\frac{n-3}{2(n-1)}J_\t} & 0 \\
	0 & 0 & \displaystyle{\frac{1}{2}L}
\end{array}\right) \, ,
$$
where $J_\theta$ and $J_\t$ are the Dirac eigenvalues of $\theta$ and 
$\tau$ on $S^{n-1}$, 
respectively 
and $L$ is the Rarita-Schwinger eigenvalue of [$\eta$] on $S^{n-1}$.\\ 
The spectrum generating relation takes the following form:\\
$$
[N,\varpi]=2\left(\nd_T+\frac{n}2\varpi\right)\,,
$$
where $\nd^{*,{\rm R}}\nd:=N$ is the Riemannian Bochner Laplacian. 
Therefore the relation (\ref{rel-1}) becomes
\begin{equation}\label{rel}
B\left(\frac{1}{2}[N,\vp]-r\vp-\X\eye[\pt,[\cD,\vp]]\right)
=\left(\frac{1}{2}[N,\vp]+r\vp+\X\eye[\pt,[\cD,\vp]]\right)B\, .
\end{equation}
As explained in detail in (\cite{Branson:97}), the recursive 
numerical spectral data come from the compressed relation of 
the above. 
%%
%%%%%%%%%%%%%%%%%%%%%%%%%%%%%%%%%%%%%%%%%%%%%%%%%%%%%%%%%%%%%%%%%%%%%%%%%
\section{Projections into isotypic summands}
%%%%%%%%%%%%%%%%%%%%%%%%%%%%%%%%%%%%%%%%%%%%%%%%%%%%%%%%%%%%%%%%%%%%%%%%%
Let us denote the $K=\mbox{Spin}(2)\times\mbox{Spin}(n)$-type 
with highest weight 
$$
(f)\otimes
(j,\tf+q\,\tf,\ldots,\tf\,,\frac{\e}{2})\, ,
$$
where $j\in\tf+q+\bbN$, 
$\e=\pm 1$, and $q=0$ or $1\,$,
by 
$$
\cV_\Xi(f,j,\tf+q\,\tf,\ldots,\tf\,,\frac{\e}{2})\, .
$$
An $\frak{s}$-map from such a $K$-type lands in the direct sum of
neighboring $K$-types (\cite{Branson:87}).\\
Consider a $\X$ spinor $\HTmx{\f}{\y}$. 
Since $\f=\Xi\y/\eye$, we have
$$
\a^0\twostac{\bullet}{\y}=\twostac{\bullet}{\Xi\y/\eye}\, .
$$
Here $\bullet$ denotes a top entry that is computable from the bottom
entry, but whose value is not needed at the moment.\\
In addition,
$$
\begin{array}{rl}
\sin t\twostac{\bullet}{\y}&=\twostac{\bullet}{\sin t\y}
=\twostac{\bullet}{-[\ptl_t,\cos t]\y}\,, \\
{\rm Proj}_{f'}\sin t\twostac{\bullet}{\y}
&=\twostac{\bullet}{\frac{f'-f}{\eye}\cos t|^{f'}_f\y}\,, \\
\sin\r\a^1\twostac{\bullet}{\y}&=\twostac{\bullet}{-\sin\r\g^1\y}
=\twostac{\bullet}{[D,\cos\r]\y}\,, \\
{\rm Proj}_b
\sin\r\a^1\twostac{\bullet}{\y}&=\twostac{\bullet}{-{\rm Proj}_b\sin\r\g^1\y}
=\twostac{\bullet}{D|^b_a\cos\r\,\y}\,,
\end{array}
$$
where $D=\g^i\nd_i$ is the Dirac operator on $S^{n-1}$.
Here $a$ and $b$ (resp., $f$ and $f'$) are abbreviated labels 
for the Spin$(n)$-types (resp., Spin$(2)$-types) in
question.\\
Note also that the compressed relations of $\vp$ between Clifford 
range part, twistor range part, and divergence part look (\ref{decom}):
\begin{equation}\label{w}
\begin{array}{l}
\vp\HTTmx{\langle \theta \rangle}{0}{0}=\HTTmx{\langle|\vp|\theta\rangle}
{0}{0}
\stackrel{\mbox{Proj}}{\longmapsto} :\HTTmx{\langle\tt\rangle}{0}{0}\, ,\\
\vp\HTTmx{0}{\{\tau\}}{0}=\HTTmx{0}{|\vp|\{\tau\}}{|\vp|\{\tau\}}
=\HTTmx{0}{C\{|\vp|\tau\}}{|\vp|\{\tau\}}
\stackrel{\mbox{Proj}}{\longmapsto}:\HTTmx{0}{C\{\tta\}}{[\eta]}\, ,\\
\vp\HTTmx{0}{0}{[\eta]}=\HTTmx{0}{|\vp|[\eta]}{|\vp|[\eta]}
\stackrel{\mbox{Proj}}{\longmapsto}:\HTTmx{0}{\{\bar{\tau}\}}
{[\tilde{\eta}]}\, ,
\end{array}
\end{equation}
where $C$ is a quantity we will compute in the following lemma.\\ 
Note that $\vp\langle \theta \rangle$ has only Clifford range pieces, 
since it is made of a spinor and fundamental tensor-spinor on $S^{n-1}$. 
On the other hand, $\vp\{\tau\}$ and $\vp[\eta]$ have no Clifford range 
pieces, since they are made of twistors on $S^{n-1}$ 
(See \cite{Branson:96, Branson:97}).
\begin{lemma} Let $\a=\cV_\X(f;j,\frac{1}{2},\cdots,\frac{1}{2},
\frac{\e}{2})$ and
$\b=\cV_\X(f';j',\frac{1}{2},\cdots,\frac{1}{2},\frac{\e'}{2})$, 
$\e=\pm 1$. Then we have
$$
|_\b\vp|_\a\{\tau\}=C_{ba}\{|_\b\vp|_\a\tau\}\, ,
$$
where
$$
C_{ba}=\frac1{\l_b(T^*T)}\left(\frac{1}{2}J_b^2+
\frac{1}{2}J_a^2-\frac{J_bJ_a}{n-1}-\frac{n(n-1)}4\right)\, ,
$$
$J_a$ (resp., $J_b$) is the Dirac eigenvalue on 
$\cV_\X(j,\frac{1}{2},\cdots,\frac{1}{2},\frac{\e}{2})$ 
(resp., $\cV_\X(j',\frac{1}{2},\cdots,\frac{1}{2},\frac{\e'}{2})$), and 
$\l_b(T^*T)$ is the eigenvalue of $T^*T$ on 
$\cV_\X(j',\frac{1}{2},\cdots,\frac{1}{2},\frac{\e'}{2})$\, over $S^{n-1}$.
\end{lemma}
\begin{proof}It suffices to show for the twistor operator $T$ on $S^{n-1}$ 
 and $\w=\cos \r$ that
$$
|_b\w|_aT\t=C_{ba}\cdot T(|_b\w|_a\t)\, .
$$
Let $D$ be the Dirac operator on $S^{n-1}$. Then
$$
\begin{array}{rl}
[D^2,\w]\t & =[\N^*\N,\w]\t \; \; 
{\rm{\textstyle  by\; Bochner\; identity}} \\
&=(\N^*\N\w)\t-2\N^k\w\N_k\t=(n-1)\w\t+2\sin\r\N_1\t\, , \\
\end{array}
$$
Also
$$
\begin{array}{rl}
T^*(\w T\t)&=-\N^j(\w\N_j\t+\tfrac1{n-1}\w\g_j D\t) \\
&=\sin\r\N_1\t+\w\N^*\N\t
+\tfrac1{n-1}\sin\r\g_1 D\t-\tfrac1{n-1}\w D^2\t \\
&=\tfrac12\left([D^2,\w]-(n-1)\w\right)\t
+\w\left(D^2-\tfrac{(n-1)(n-2)}4
\right)\t+\tfrac1{n-1}[\w,D]D\t \\
&\qquad -\tfrac1{n-1}\w D^2\t \quad 
{\rm{\textstyle by\; the\; above\; and\; Bochner\;
identity}} \\
&=\tfrac12 D^2(\w\t)+\tfrac12 \w D^2\t
-\tfrac1{n-1}D(\w D\t)-\tfrac{n(n-1)}4\w\t.  
\end{array}
$$
Therefore 
$$
\begin{array}{ll}
|_b\w|_aT\t&=T\left(\dfrac{1}{\l_b(T^*T)}T^*(|_b\w|_aT\t)\right)\\
    &=T\left(\dfrac{1}{\l_b(T^*T)}\left(\frac{1}{2}J_b^2+\frac{1}{2}J_a^2
    -\frac{1}{n-1}J_bJ_a-\frac{n(n-1)}{4}\right)|_b\w|_a\t\right)\, .
\end{array}
$$
\end{proof}
\begin{rmk}Eigenvalues of $D$ and $T^*T$ on $S^{n-1}$ are known 
due to Branson (\cite{Branson:99}). 
\end{rmk}
\noindent With the above (\ref{w}) at hand, we get
\begin{equation}\label{Dw}
\begin{array}{l}
|_\b[\cD,\vp]|_\a\langle \theta \rangle=\HTTmx{(\cD_{11}^\b-\cD_{11}^\a)
\langle\tt\rangle}
{(\cD_{21}^\b-C_{ba}\cD_{21}^\a)\{\tt\}}{-\cD_{21}^\a[\eta]}, 
\mbox{ where }
\left\{\htmx{\langle \tt \rangle=|_\b\vp|_\a\langle\theta\rangle}
{[\eta]=|_\b\vp|_\a\{\theta\}}\right. ,\\
|_\b[\cD,\vp]|_\a\{\tau\}=
\HTTmx{(C_{ba}\cD_{12}^\b-\cD_{12}^\a)\langle\tta\rangle}
{C_{ba}(\cD_{22}^\b-\cD_{22}^\a)\{\tta\}}
{(\cD_{33}^\b-\cD_{22}^\a)[\eta]}, 
\mbox{ where }
\left\{\htmx{\{\tta\}=|_\b\vp|_\a\{\tau\}}{[\eta]=|_\b\vp|_\a\{\tau\}}
\right. ,\mbox{ and}\\
|_\b[\cD,\vp]|_\a[\eta]=\HTTmx{\cD_{12}^\b\langle\bar{\tau}\rangle}
{(\cD_{22}^\b-\cD_{33}^\a)\{\bar{\tau}\}}
{(\cD_{33}^\b-\cD_{33}^\a)[\tilde{\eta}]}
, \mbox{ where } 
\left\{\htmx{\{\bar{\tau}\}=|_\b\vp|_\a[\eta]}
{[\tilde{\eta}]=|_\b\vp|_\a[\eta]}\right. .
\end{array}
\end{equation}
Here we use subscripts to refer to the specific entries of the $\cD$ 
and superscripts to indicate where these entries are computed.\\
Let us now consider the compressed relation of (\ref{rel}) 
between neighboring $K$-types.\\ \\
%%%%%%%%%%%%%%%%%%%%%%%%%%%%%%%%%%%%%%%%%%%%%%%%%%%%%%%%%%%%%%%%%%%
{\bf Case 1: Multiplicity 2 $\leftrightarrow$ 1}
%%%%%%%%%%%%%%%%%%%%%%%%%%%%%%%%%%%%%%%%%%%%%%%%%%%%%%%%%%%%%%%%%%%
$$
\a=\cV_\X(f;j,\frac{1}{2},\cdots,\frac{1}{2},\frac{\e}{2}) 
\leftrightarrow
\b=\cV_\X(f';j,\frac{3}{2},\frac{1}{2},\cdots,
\frac{1}{2},\frac{\e}{2})\, .
$$
Note that the operator $B$ in block form looks
$$
B=\left(\begin{array}{ccc}B_{11}&B_{12}&0\\
                    B_{21}&B_{22}&0\\
		    0&0&B_{33}
  \end{array}\right)\, .
$$
With
$$
|_{\a}N|_{\b}=f^2-f'^2-(n-2)\mbox{ and }
|_{\b}N|_{\a}=-|_{\a}N|_{\b}
$$
and (\ref{Dw}), we get $\a\rightarrow\b$ transition quantities
$$
\begin{array}{cc}
\b\rightarrow\a :&
\HHmx{B_{11}^{\a}}{B_{12}^{\a}}{B_{21}^{\a}}{B_{22}^{\a}}
\HTmx{A_1}{E^-}
=B_{33}^{\b}\HTmx{-A_1}{E^+}\, \mbox{ and}\\
\a\rightarrow\b :&
\THmx{A_2}{-E^-}\HHmx{B_{11}^{\a}}
{B_{12}^{\a}}{B_{21}^{\a}}{B_{22}^{\a}}
=B_{33}^{\b}\THmx{-A_2}{-E^+}\, ,
\end{array}
$$
where
$$
\begin{array}{l}
A_1:=\X(f-f')\cD_{12}^\a\, ,\\
A_2:=-\X(f-f')\cD_{21}^\a\, ,\\
E^-:=\frac{1}{2}(f^2-f'^2)-\frac{n-2}{2}-r+\X(f-f')
(\cD_{22}^\a-\cD_{33}^\b)\, 
,\\
E^+:=\frac{1}{2}(f^2-f'^2)-\frac{n-2}{2}+r-\X(f-f')
(\cD_{22}^\a-\cD_{33}^\b)\, .
\end{array}
$$
In particular,  we can write all $2\times 2$ entries of $B^{\a}$ 
in terms of $B_{21}^{\a}$ and $B_{33}^{\b}$:
\begin{equation}\label{m1andm2}
\begin{array}{l}
B_{11}^{\a}=(E^-B_{21}^{\a}-A_2B_{33}^{\b})
/{A_2}\, ,\\
B_{12}^{\a}=-A_1B_{21}^{\a}/A_2\, ,
\mbox{ and}\\
B_{22}^{\a}=(-A_1B_{21}^{\a}+E^+B_{33}^{\b})
/E^-\, .
\end{array}
\end{equation}
Thus if we can express $B_{21}^{\a}$ in terms of $B_{33}^{\b}$, 
we can completely determine all entries in the $2\times2$ block.
\\ \\
%%%%%%%%%%%%%%%%%%%%%%%%%%%%%%%%%%%%%%%%%%%%%%%%%%%%
\noindent {\bf Case 2: Multiplicity 2 $\leftrightarrow$ 2}
%%%%%%%%%%%%%%%%%%%%%%%%%%%%%%%%%%%%%%%%%%%%%%%%%%%%
$$
\a=\cV_\X(f;j,\frac{1}{2},\cdots,\frac{1}{2},\frac{\e}{2}) 
\rightarrow
\b=\cV_\X(f';j'\frac{1}{2},\cdots,\frac{1}{2},\frac{\e'}{2})\, .
$$
Here we have
$$
|_{\b}N|_{\a}=f'^2-f^2+J_b^2-J_a^2\, .
$$
So using (\ref{Dw}), we get the transition quantities
\begin{equation}\label{m2andm2}
\HHmx{B_{11}^{\b}}{B_{12}^{\b}}{B_{21}^{\b}}{B_{22}^{\b}}
\HHmx{F_1^-}{G_2}{G_1}{C_{ba}F_2^-}
=\HHmx{F_1^+}{-G_2}{-G_1}{C_{ba}F_2^+}
\HHmx{B_{11}^{\a}}{B_{12}^{\a}}{B_{21}^{\a}}{B_{22}^{\a}}\, ,
\end{equation}
where
$$
\begin{array}{l}
F_1^-:=\frac{1}{2}(f'^2-f^2)+\frac{1}{2}(J_b^2-J_a^2)
-r+\X(f'-f)(\cD_{11}^\b-\cD_{11}^\a)\, ,\\
F_1^+:=\frac{1}{2}(f'^2-f^2)+\frac{1}{2}(J_b^2-J_a^2)
+r-\X(f'-f)(\cD_{11}^\b-\cD_{11}^\a)\, ,\\
F_2^-:=\frac{1}{2}(f'^2-f^2)+\frac{1}{2}(J_b^2-J_a^2)
-r+\X(f'-f)(\cD_{22}^\b-\cD_{22}^\a)\, ,\\
F_2^+:=\frac{1}{2}(f'^2-f^2)+\frac{1}{2}(J_b^2-J_a^2)
+r-\X(f'-f)(\cD_{22}^\b-\cD_{22}^\a)\, ,\\
G_1:=\X(f'-f)(\cD_{21}^\b-C_{ba}\cD_{21}^\a)\, , 
\mbox{ and}\\
G_2:=\X(f'-f)(C_{ba}\cD_{12}^\b-\cD_{12}^\a)\, .
\end{array}
$$
Therefore we get determinant quotients of $B$ 
on multiplicity 2 part.\\
Note the following diagram of neighboring multiplicity 2 
isotypic summands centered at  
$\cV_\X(f;j,\frac{1}{2},\cdots,\frac{1}{2},\frac{\e}{2})$:
$$
\begin{array}{lcccl}
\cV_\X(f-1;j+1,\frac{1}{2},\cdots,\frac{1}{2},\frac{\e}{2})
&&&&
\cV_\X(f+1;j+1,\frac{1}{2},\cdots,\frac{1}{2},\frac{\e}{2})\\
&\nwarrow&&\nearrow&\\
\cV_\X(f-1;j,\frac{1}{2},\cdots,\frac{1}{2},-\frac{\e}{2})
&\leftarrow&\bullet&\rightarrow&
\cV_\X(f+1;j,\frac{1}{2},\cdots,\frac{1}{2},-\frac{\e}{2})\\
&\swarrow&&\searrow&\\
\cV_\X(f-1;j-1,\frac{1}{2},\cdots,\frac{1}{2},\frac{\e}{2})
&&&&
\cV_\X(f+1;j-1,\frac{1}{2},\cdots,\frac{1}{2},\frac{\e}{2})\, .
\end{array}
$$
The determinant quotients corresponding to the above 
diagram are:
\begin{equation}\label{det-m2}
\left(\begin{array}{ll}
\frac{\left(-f+J+1-\X+r+\frac{\e}{2}\X\right)
\left(-f+J+1+\X+r+\frac{\e}{2}\X\right)}
{\left(-f+J+1-\X-r-\frac{\e}{2}\X\right)
\left(-f+J+1+\X-r-\frac{\e}{2}\X\right)}&
\frac{\left(f+J+1-\X+r-\frac{\e}{2}\X\right)
\left(f+J+1+\X+r-\frac{\e}{2}\X\right)}
{\left(f+J+1-\X-r+\frac{\e}{2}\X\right)
\left(f+J+1+\X-r+\frac{\e}{2}\X\right)}\\&\\
\frac{\left(-f+\frac{1}{2}-\X+r-\e\X J\right)
\left(-f+\frac{1}{2}+\X+r-\e\X J\right)}
{\left(-f+\frac{1}{2}-\X-r+\e\X J\right)
\left(-f+\frac{1}{2}+\X-r+\e\X J\right)}&
\frac{\left(f+\frac{1}{2}-\X+r+\e\X J\right)
\left(f+\frac{1}{2}+\X+r+\e\X J\right)}
{\left(f+\frac{1}{2}-\X-r-\e\X J\right)
\left(f+\frac{1}{2}+\X-r-\X J\right)}\\&\\
\frac{\left(-f-J+1-\X+r-\frac{\e}{2}\X\right)
\left(-f-J+1+\X+r-\frac{\e}{2}\X\right)}
{\left(-f-J+1-\X-r+\frac{\e}{2}\X\right)
\left(-f-J+1+\X-r+\frac{\e}{2}\X\right)}&
\frac{\left(f-J+1-\X+r+\frac{\e}{2}\X\right)
\left(f-J+1+\X+r+\frac{\e}{2}\X\right)}
{\left(f-J+1-\X-r-\frac{\e}{2}\X\right)
\left(f-J+1+\X-r-\frac{\e}{2}\X\right)}
\end{array}\right)\, ,
\end{equation}
where $J=\e J_a$. 
\\
And these data can be put into the 
following gamma function expression:
$$
\frac{1}{4}\bullet
\frac{\G\left(\frac12(f+J+r-\frac{\e}{2}\Xi)\right)
\G\left(\frac12(-f+J+r+\frac{\e}{2}\Xi)\right)}
{\G\left(\frac12(f+J-r+\frac{\e}{2}\Xi)\right)
\G\left(\frac12(-f+J-r-\frac{\e}{2}\Xi)\right)}
$$
$$
\bullet
\frac{\G\left(\frac12(f+J+2+r-\frac{\e}{2}\Xi)\right)
\G\left(\frac12(-f+J+2+r+\frac{\e}{2}\Xi)\right)}
{\G\left(\frac12(f+J+2-r+\frac{\e}{2}\Xi)\right)
\G\left(\frac12(-f+J+2-r-\frac{\e}{2}\Xi)\right)}\, .
$$
\\ \\
%%%%%%%%%%%%%%%%%%%%%%%%%%%%%%%%%%%%%%%%%%%%%%%%%%%%%%
{\bf Case 3: Multiplicity 1 $\leftrightarrow$ 1}
%%%%%%%%%%%%%%%%%%%%%%%%%%%%%%%%%%%%%%%%%%%%%%%%%%%%%%
$$
\a=\cV_\X(f;j,\frac{3}{2},\frac{1}{2},\cdots,\frac{1}{2},
\frac{\e}{2}) 
\leftarrow
\b=\cV_\X(f';j'\frac{3}{2},\frac{1}{2},\cdots,\frac{1}{2},
\frac{\e'}{2})
\, .
$$
Again we have
$$
|_{\a}N|_{\b}=f^2-f'^2+J_a^2-J_b^2\, .
$$
And the transition quantities are
\begin{equation}\label{m1andm1}
B_{33}^\a P^-=P^+ B_{33}^\b\, ,
\end{equation}
where
$$
\begin{array}{l}
P^-:=\frac{1}{2}(f^2-f'^2)+\frac{1}{2}(J_a^2-J_b^2)
-r+\X(f-f')(\cD^\a_{33}-\cD^\b_{33})\mbox{ and}\\
P^+:=\frac{1}{2}(f^2-f'^2)+\frac{1}{2}(J_a^2-J_b^2)
+r-\X(f-f')(\cD^\a_{33}-\cD^\b_{33})\, .
\end{array}
$$
The diagram of neighboring multiplicity 1 isotypic 
summands centered at 
$$
\cV_\X(f;j,\frac{3}{2},\frac{1}{2},\cdots,\frac{1}{2},\frac{\e}{2})
$$ 
looks:
$$
\begin{array}{lcccl}
\cV_\X(f-1;j+1,\frac{3}{2},\frac{1}{2},\cdots,\frac{1}{2},\frac{\e}{2})
&&&&
\cV_\X(f+1;j+1,\frac{3}{2},\frac{1}{2},\cdots,\frac{1}{2},\frac{\e}{2})\\
&\nwarrow&&\nearrow&\\
\cV_\X(f-1;j,\frac{3}{2},\frac{1}{2},\cdots,\frac{1}{2},-\frac{\e}{2})
&\leftarrow&\bullet&\rightarrow&
\cV_\X(f+1;j,\frac{3}{2},\frac{1}{2},\cdots,\frac{1}{2},-\frac{\e}{2})\\
&\swarrow&&\searrow&\\
\cV_\X(f-1;j-1,\frac{3}{2},\frac{1}{2},\cdots,\frac{1}{2},\frac{\e}{2})
&&&&
\cV_\X(f+1;j-1,\frac{3}{2},\frac{1}{2},\cdots,\frac{1}{2},\frac{\e}{2})
\, .
\end{array}
$$
And the eigenvalue quotients are:
$$
\left(\begin{array}{ll}
\dfrac{-f+J+1+r+\frac{\e}{2}\X}{-f+J+1-r-\frac{\e}{2}\X}&
\dfrac{f+J+1+r-\frac{\e}{2}\X}{f+J+1-r+\frac{\e}{2}\X}\\&\\
\dfrac{-f+\frac{1}{2}+r-\e\X J}{-f+\frac{1}{2}-r+\e\X J}&
\dfrac{f+\frac{1}{2}+r+\e\X J}{f+\frac{1}{2}-r-\e\X J}\\&\\
\dfrac{-f-J+1+r-\frac{\e}{2}\X}{-f-J+1-r+\frac{\e}{2}\X}&
\dfrac{f-J+1+r+\frac{\e}{2}\X}{f-J+1-r-\frac{\e}{2}\X}
\end{array}\right)\, ,
$$
where $J=\e J_a$.\\
Thus, following the normalization on the multiplicity 2 part, 
we get the spectral function on the multiplicity 1 part:
\begin{equation}\label{SpecFcn}
Z(r;f,J,\Xi\e)=\frac{\e}{2}\Xi
\frac{\G\left(\frac12(f+J+1+r-\frac{\e}{2}\Xi)\right)
\G\left(\frac12(-f+J+1+r+\frac{\e}{2}\Xi)\right)}
{\G\left(\frac12(f+J+1-r+\frac{\e}{2}\Xi)\right)
\G\left(\frac12(-f+J+1-r-\frac{\e}{2}\Xi)\right)}\,.
\end{equation}
In particular,
$$
Z(\frac{1}{2},f,J,\X\e)=-\frac{1}{4}(f-\X\e J)
=\frac{1}{4}\eye\mbox{ eig}(E\cR;f,J,\X\e)\, ,
$$
where $E\cR$ is the exchanged Rarita-Schwinger operator.
%%%%%%%%%%%%%%%%%%%%%%%%%%%%%%%%%%%%%%%%%%%%%%%%%%%%%%%%%%%%%%%%%%%%%%%%
\section{Interface between multiplicity 1 and 2 parts}
%%%%%%%%%%%%%%%%%%%%%%%%%%%%%%%%%%%%%%%%%%%%%%%%%%%%%%%%%%%%%%%%%%%%%%%%
Consider the following diagram:
$$
\begin{array}{ccc}
\a_1=\cV_\X(f;j,\frac{1}{2},\cdots,\frac{1}{2},\frac{\e}{2}) 
&\rightarrow&
\a_2=\cV_\X(f+1;j+1,\frac{1}{2},\cdots,\frac{1}{2},\frac{\e}{2})\\
&&\\
\updownarrow & &\updownarrow \\
&&\\
\b_1=\cV_\X(f+1;j,\frac{3}{2},\frac{1}{2},\cdots,\frac{1}{2},\frac{\e}{2}) 
&\leftarrow&
\b_2=\cV_\X(f;j+1,\frac{3}{2},\frac{1}{2},\cdots,\frac{1}{2},\frac{\e}{2})\, .
\end{array}
$$
Then (\ref{m2andm2}) reads
$$
B^{\a_2}M_1=M_2B^{\a_1}\, .
$$
So
$$
{\rm det}B^{\a_2}=\dfrac{{\rm det}M_2}{{\rm det}M_1}{\rm det}B^{\a_1}\, .
$$
Note that $\dfrac{{\rm det}M_2}{{\rm det}M_1}$ is a determinant 
quotient computed in (\ref{det-m2}).\\
From (\ref{m1andm2}), we get a relation between $B_{12}$ and $B_{33}$:
$$
\begin{array}{ll}
\mbox{det}\HHmx{B_{11}}{B_{12}}{B_{21}}{B_{22}}&
=B_{11}B_{22}-B_{12}B_{22}\\
&=-\dfrac{1}{A_2E^-}B_{33}\left(B_{33}A_2E^+-(E^-E^++A_1A_2)B_{21}
\right)\, .
\end{array}
$$
We can also compare $(2,1)$ entries of both sides in (\ref{m2andm2}). 
Applying (\ref{m1andm2}) and (\ref{m1andm1}) to the both relations, 
we can finally write $B_{21}$ in terms of $B_{33}$
with a ``big'' help 
from computer algebra package.
\par
$2\times 2$ block on 
$$
\cV_\X(f;j,\frac{1}{2},\cdots,\frac{1}{2},\frac{\e}{2})
$$
in terms of $(3,3)$ 
$$
\cV_\X(f+1;j,\frac{3}{2},\frac{1}{2},\cdots,
\frac{1}{2},\frac{\e}{2})
$$
is:
\begin{equation}\label{block}
\HHmx{\dfrac{4C_1C_2}{(n-1)C_3C_4}-1}
{\dfrac{-2(n-2)\Xi C_5C_2}{(n-1)^2C_3C_4}}
{\dfrac{8n\Xi C_2}{C_3C_4}}
{\dfrac{-4C_5C_2}{(n-1)C_1C_3C_4}+\dfrac{C_6}{C_1}}
\bullet
Z(r;f+1,J,\X\e)\, ,
\end{equation}
where
$$
\begin{array}{l}
C_1=2fn-2f-2n+1+n^2+2rn-2r-2\Xi J_a \, , \\
C_2=2fr+\Xi J_a  \, ,\\
C_3=n-1+2r \, , \\
C_4=(2f+2r-\Xi+2J_a)(2f+2r+\Xi-2J_a)\, ,  \\
C_5=(n-1+2J_a)(n-1-2J_a) \, ,\mbox{ and} \\
C_6=2fn-2f-2n+1+n^2-2rn+2r+2\Xi J_a\, .  
\end{array}
$$
\begin{rmk}In particular, if $r=\frac{1}{2}$ 
and $(3,3)$ entry 
$$
\eye f-\eye\X\e J
$$
of the exchanged Rarita-Schwinger operator is put 
into the above formula, 
we recover the other  
$2\times 2$ entries
$$
\left(\begin{array}{cc}
-\dfrac{n-2}{n}\eye\left(f+\dfrac{n+1}{n-1}\X\e J\right)
&-\dfrac{2\eye\X}{n(n-1)}\left(\dfrac{(n-1)(n-2)}{4}
-\dfrac{n-2}{n-1}J^2\right)\\
2\eye\X &\eye f-\dfrac{n-3}{n-1}\eye\X\e J
\end{array}\right).
$$
\end{rmk}
%%%%%%%%%%%%%%%%%%%%%%%%%%%%%%%%%%%%%%%%%%%%%%%%%%%%%%%%%%%%%%%%%%%%%%%%%%
%%%%%%%%%%%%%%%%%%%%%%%%%%%%%%%%%%%%%%%%%%%%%%%%%%%%%%%%%%%%%%%%%%%%%%%%%%

\end{document}